\documentclass[11pt, draft]{paper}        

\usepackage{amsfonts}
\usepackage{amsmath}
\usepackage{amsthm}
\usepackage{amscd}
\usepackage{euscript}
\usepackage{amssymb}

\usepackage{psfrag}
\usepackage{graphicx}   
\usepackage{rotating}
\usepackage{color}
\usepackage{mathrsfs}
\usepackage{xy}
\xyoption{all}

\theoremstyle{plain}
\newtheorem{thm}[equation]{Theorem}

\newtheorem{lem}[equation]{Lemma}
\newtheorem{prop}[equation]{Proposition}
\newtheorem{conj}[equation]{Conjecture}
\theoremstyle{definition}
\newtheorem{defn}[equation]{Definition}
\theoremstyle{remark}

\numberwithin{equation}{section}

\newenvironment{prf}[1][]
  {\medskip\par\noindent{\bf Proof#1. }}
  {\nopagebreak\par\rightline{$\Box$} \medskip}

\newcommand{\CC}{\mathbb C}

\newcommand{\GG}{\mathbb G}

\newcommand{\PP}{{\mathbb P}}

\newcommand{\LG}{\mathbb{LG}}

 \DeclareMathOperator{\rk}{rk}

\DeclareMathOperator{\Sp}{Sp}
\DeclareMathOperator{\Tang}{T}

\newcommand{\SymplFormMtrx}[1]{
        \begin{bmatrix}
                   0 & I_{#1}\\
                -I_{#1} &   0
        \end{bmatrix}}

\newcommand{\pder}[1]{\frac{\partial}{\partial #1}}

\newcommand{\set}[1]{\left\{#1\right\}}
\newcommand{\fromto}[2]{#1, \dotsc, #2}
\newcommand{\setfromto}[2]{\set{\fromto{#1}{#2}}}

\newcommand{\MaxVn}{8}

\def\Wedge#1{{\textstyle{\bigwedge\nolimits}^{\! #1}}}

\newcommand{\RotateUp}[1]{\rotatebox{90}{\scriptsize\ensuremath{#1}}}

\title{Secants of Lagrangian Grassmannians}

\author{Ada Boralevi
\thanks{A.~Boralevi's e-mail address: boralevi@math.tamu.edu} 
\and
Jaros\l{}aw Buczy\'{n}ski
\thanks{J.~Buczynski is supported by Marie Curie International Outgoing Fellowship;
 e-mail address: jabu@mimuw.edu.pl}}
\institution{Texas A\&M University, College Station, TX, USA}

\begin{document}

\maketitle
\abstract{%
We study the dimensions of secant varieties of Grassmannian
of Lagrangian subspaces in a symplectic vector space.
We calculate these dimensions for third and fourth secant varieties.
Our result is obtained by providing a normal form for four general points on such a Grassmannian
and by explicitly calculating the tangent spaces at these four points.}

\smalltableofcontents


\section{Introduction}


Let $X \subset \PP^N$ be a non-degenerate projective variety. The $r$-th secant variety $\sigma_r(X)$ is 
defined to be the closure of the union of linear spans of all the $r$-tuples of points lying on $X$.

It is a long standing and well established problem to calculate properties of secants of certain varieties, 
in particular homogeneous spaces in their homogeneous embeddings. 
The dimension is the simplest among those investigated properties,
yet even for the easiest homogeneous spaces,
the calculation of dimension is a highly non-trivial problem, and there is an extensive related literature. 
A well known classification of defective secants to Veronese embeddings of $\PP^n$ was completed 
in a series of papers by Alexander and Hirschowitz \cite{alexander_hirschowitz_polynomial_interpolation}. 
There are corresponding conjecturally complete lists of defective secants to Segre products $\PP^{n_1} \times \dotsm \times \PP^{n_k}$ \cite{abo_ottaviani_peterson_segre} 
and to ordinary Grassmannians $\GG(k,n)$ (see \cite{abo_ottaviani_peterson_grassmannian_of_planes}, \cite{catalisano_geramita_gimigliano} and 
\cite{baur_draisma_degraaf}), whereas for Segre-Veronese varieties even such a conjectural classification is missing (see \cite{abo_brambilla} and numerous references therein).

In this paper we undertake the study of dimensions of secant varieties of Lagrangian Grassmannians $\LG(n,2n)$
in their minimal homogeneous embeddings.
These are projective varieties parametrising dimension $n$ isotropic subspaces of a symplectic vector space $V$ 
of dimension $2n$.

\begin{thm}\label{thm third and fourth secant}
Suppose $n \ge 4$, and $r =3$ or $r=4$.
Then:
\begin{itemize}
 \item If $n=4$, $r=3$, then $\dim \sigma_3(\LG(4,8)) = 31 = (3*11 -1) -1$.
 \item If $n=4$, $r=4$, then $\dim \sigma_4(\LG(4,8)) = 39 = (4*11 -1) -4$. 
 \item If $n\ge 5$, then $\sigma_3(\LG(n,2n))$  and $\sigma_4(\LG(n,2n))$ always have the expected dimension,
       namely $r(d+1) -1$, where $d= \dim \LG(n, 2n)= \binom{n+1}{2}$
\end{itemize}
\end{thm}
The cases $n \le 3$ or $r=2$ are also explained in our paper, but these were known before, see Section~\ref{known cases}.

The proof of the theorem is split into the case $r=3$ and $r=4$ and explained in
 Sections~\ref{section third secant case} and \ref{section fourth secant case}, respectively.
The idea  is to calculate a normal form for four general points of a Lagrangian Grassmannian 
(see Proposition~\ref{prop normal form}) 
and then perform an explicit calculation of generators of the affine tangent spaces at this normalised general points.
By application of Terracini Lemma (see Lemma~\ref{Terracini} below)
the dimension of the secant variety is determined by the dimension of the sum of those affine tangent spaces.
Thus the theorem boils down to the calculation of a rank of a certain matrix,
whose rows are the generators of the tangent spaces.

In \S\ref{section_computations} we also present the results of some computational experiments in Magma \cite{magma}
in small dimensions. As a conclusion we dare to conjecture:
\begin{conj}\label{conjecture}
  The secant variety $\sigma_r(\LG(n,2n))$ has the expected dimension except for the case $n=4$ and $r=3$
  (defect $1$) and $r=4$ (defect $4$).
\end{conj}
As a corollary of Theorem~\ref{thm third and fourth secant} and Proposition~\ref{prop small values of n} we know 
the conjecture is true both for $r \le 4$ and also for $n \le \MaxVn$. 

We point out that Conjecture~\ref{conjecture} is somehow ``expected'' from the related conjectural classification of defective secants to ordinary Grassmannians. 
We quote from \cite[Conjecture 4.1]{baur_draisma_degraaf}, yet such a conjecture was long before believed to be true.

\begin{conj}
  Let $k \geq 3$. Then $\sigma_r(\GG(k,n))$ has the expected dimension except for the cases $\sigma_3(\GG(3,7))$ (defect $1$), 
$\sigma_3(\GG(4,8))$ and $\sigma_4(\GG(4,8))$ (defect $1$ and $4$ respectively) and $\sigma_4(\GG(3,9))$ (defect $2$).
 \end{conj}

\subsection*{Acknowledgements}

We would like to thank Giorgio Ottaviani for suggesting this topic and his precious comments.


\section{Known cases in small dimension}\label{known cases}


We briefly review what is known in low dimension.

The first Lagrangian Grassmannian $\LG(2,4)$ has dimension $d=3$ and is a quadric hypersurface in $\PP^4$, 
so $\sigma_2(\LG(2,4))$ fills the ambient $\PP^4$.

The variety $\LG(3,6)$ has dimension $d=6$ and 
\[
  \LG(3,6) \subset \PP^{13} = \PP\left(\bigoplus_{i=0}^3 (S^2 \Wedge{i} \CC^3)\right).
\]
Again the secant variety $\sigma_2(\LG(3,6))=\PP^{13}$ fills the ambient space,
as follows from Lemma~\ref{tg space p1 and p2} below.
This statement was known before ---
for instance it is contained in the proof of \cite[Prop.~17(2)]{landsbergmanivel04},
since $\LG(3,6)$ is a Legendrian subvariety in $\PP^{13}$. 

In general, $\sigma_2 (\LG(n,2n))$ has the expected dimension for all $n$, as observed in Lemma~\ref{tg space p1 and p2}. 
This must have been known before, although we are unable to find an explicit reference in the literature.


\section{Notations and definitions}


Throughout the paper we work over an algebraically closed base field of characteristics zero. 

If $p_1,\ldots,p_r$ are points in $\PP^N$ we let $\langle p_1,\ldots,p_r\rangle$ denote their linear span.
If $X \subset \PP^N$, then by definition the $r$-th secant variety $\sigma_r(X)$ is: 
\[
 \sigma_r(X)=\overline{\bigcup_{p_1,\ldots,p_r \in X} \langle p_1,\ldots,p_r \rangle}.
\]
If $X \subset \PP^N$ is  non-degenerate and $\dim X = d$,
 then the dimension of $\sigma_r(X)$ cannot exceed $\min\{N,r(d+1)-1\}$.

We use a standard terminology, which is summarized in the following definition:

\begin{defn} Let $X \subset \PP^N$ be a non-degenerate variety of dimension $d$.
  \begin{enumerate}
    \item If $\dim \sigma_r(X)=\min\{N,r(d+1)-1\}$ 
          we say that $\sigma_r(X)$ \emph{has the expected dimension}.
    \item If $\dim \sigma_r(X) < \min\{N,r(d+1)-1\}$ 
          we say that $X$ is \emph{$r$-defective}, or that it \emph{has a defective $r$-th secant variety}.
    \item If $X$ is $r$-defective, its \emph{defect} is the difference $r(d+1)-1 - \dim \sigma_r(X)$.
  \end{enumerate}
\end{defn}

The main tool used to compute the dimension of secant varieties is the following well-known lemma by Terracini,
 see \cite[Proposition 1.10]{zak_tangents}:

\begin{lem}[Terracini Lemma]\label{Terracini}
Let $p_1,\ldots,p_r$ be general points in $X$ and let $z$ be a general point of $ \langle p_1,\ldots,p_r \rangle$.
 Then the affine tangent space to $\sigma_r(X)$ at $z$ is given by 
\[
  \hat{\Tang}_z\sigma_r(X)=  \hat{\Tang}_{p_1}X  + \dotsb + \hat{\Tang}_{p_r}X 
\]
where $\hat{\Tang}_{p_i}X$ denotes the affine tangent space to $X$ at $p_i$.
\end{lem}

We call $\LG(n,2n)$ the Lagrangian Grassmannian of dimension $n$ 
Lagrangian subspaces of a complex symplectic vector space $V$ of dimension $2n$. 
It may be identified with the homogeneous space $\Sp(2n)/P(\alpha_n)$ of dimension $d=\binom{n+1}{2}$. Here 
$\alpha_n$ is the last simple root of the Lie algebra of $\Sp(2n)$, and $P(\alpha_n) \leq \Sp(2n)$ is 
the parabolic subgroup obtained by removing the root $\alpha_n$. 
(We use the ordering of the roots as in \cite{bourbakitables}.)

We fix a symplectic basis of $V$, by which we mean that the matrix of the symplectic form in this basis is 
   $J=\SymplFormMtrx{n}$.
An element $p \in \LG(n,2n)$ is a vector subspace of $V$, and it can be represented by 
its basis identified with a $2n \times n$ matrix of rank $n$ that looks like this:
\[
  B=\begin{bmatrix}
       B_1\\B_2
    \end{bmatrix}
\]
with  $n \times n$ square matrices $B_1$ and $B_2$.
We will say $p = \theta(B)$.
Note $p = \theta(B)$ is an isotropic subspace if and only if $B^T J B=0$,
that is if and only if $B_1^T B_2 = B_2^T B_1$.
 
A presentation of a vector space with a basis clearly requires a choice,
so there is a $GL(n)$-action on all possible choices of $B$
that define the same $p\in \LG(n,2n)$:
\[
 B \cdot g = 
    \begin{bmatrix}
       B_1 g \\B_2 g
    \end{bmatrix}
\]

If we restrict our attention to an open affine neighborhood of 
$
	\begin{bmatrix}
                I_n\\
		0
	\end{bmatrix}
$,
consisting of those
$
  B=\begin{bmatrix}
       B_1\\B_2
    \end{bmatrix}
$
for which $B_1$ is invertible,
then every point of $\LG(n,2n)$
in this neighborhood is represented uniquely by a matrix of the form 
$
	\begin{bmatrix}
                I_n\\
		A
	\end{bmatrix}
$
where $A$ is a symmetric $n \times n$ matrix and $A = B_2 {B_1}^{-1}$.
Throughout the paper we use the convention that $A$ has the form
$$
\left[
\begin{array}{cccc}
2a_{11}&a_{12}&\ldots&\\
a_{12}&2a_{22}&&\\
&&\ddots&\\
&&&2a_{nn}
\end{array}
\right].
$$
Here $a_{ij}$ are going to be treated as local coordinates on $\LG(n, 2n)$ around 
$\theta\left(
	\begin{bmatrix}
                I_n\\
		0
	\end{bmatrix}
 \right)$.


\section{Normal forms}\label{normal forms}


The Lagrangian Grassmannian $\LG(n, 2n)$ is a homogeneous space of dimension 
$d = \binom{n+1}{2}$ with the transitive action of $\Sp(2n)$, a group of dimension 
$n(2n+1) = 4d - n$. 
Thus we expect that the quadruples of general points of $\LG(n, 2n)$ up to the action of $\Sp(2n)$ are parametrised by an $n$-dimensional family.
This is the case and we explicitly describe this family in the following proposition.

\begin{prop}\label{prop normal form}
   Let $p_1, p_2, p_3, p_4 \in \LG(n,2n)$ be four general points.
   Then there exists a choice of symplectic coordinates on $V$,
   such that:
   \begin{align*}
       p_1 & = 
	\theta\left(
	  \begin{bmatrix}
                I_n\\
		0
	  \end{bmatrix} \right) &
       p_2 & = 
	\theta\left(
	  \begin{bmatrix}
                0\\
		I_n
	  \end{bmatrix} \right) &
       p_3 & = 
	\theta\left(
	  \begin{bmatrix}
                I_n\\
		I_n
	  \end{bmatrix}  \right)&
       p_4 & = 
	\theta\left(
	  \begin{bmatrix}
                I_n\\
		Q_n
	  \end{bmatrix}  \right)
   \end{align*}
   where $Q_n=diag(q_1,\ldots,q_n)$ is a general $n \times n$ diagonal matrix.
\end{prop}

\begin{prf}
   By homogeneity of $\LG(n,2n)$, the choice of the first point is arbitrary.
   General Lagrangian subspaces are pairwise disjoint.
   For two disjoint Lagrangian subspaces $p_1$, $p_2$, 
   their direct sum is $p_1\oplus p_2 = V$ by dimension count.
   Also the symplectic form identifies $p_2$ with ${p_1}^*$ in such a way 
   that the symplectic form on $ V$ is the standard symplectic form on  $p_1 \oplus {p_1}^*$.
   Choose any basis of $p_1$ and the dual basis of ${p_1}^* \simeq p_2$ and this gives the 
   normal form of $p_2$.
 
   Having fixed $p_1$ and $p_2$, 
   we still have a large subgroup of $Sp(2n)$ acting on $\LG(n,2n)$ and 
   preserving $p_1$ and $p_2$.
   Namely, this is $GL(n)$ acting as follows.  
   For $g\in GL(n)$, the following matrix 
   $
    \begin{bmatrix}
        g^{-1} &    0\\
        0      & g^T
    \end{bmatrix}
   $
   is the symplectomorphism representing this action.
   If $p_3 \in \LG(n,2n)$ is a general element, then we may assume it is in the open affine neighborhood 
   of $p_1$ and thus it is of the form 
   $
	\theta\left(\begin{bmatrix}
                I_n\\
		A
	\end{bmatrix} \right)
   $ for some symmetric matrix $A$. 
   By generality we may also assume $A$ is non-degenerate.
   Now:
   \[
     \theta\left(\begin{bmatrix}
        g^{-1} &    0\\
        0      & g^T
    \end{bmatrix}
    \begin{bmatrix}
                I_n\\
                A
    \end{bmatrix}\right)
    =
    \theta\left(\begin{bmatrix}
                g^{-1}\\
                 g^T A
    \end{bmatrix}\right)
    =
    \theta\left(\begin{bmatrix}
                I_n \\
                g^T A g
    \end{bmatrix}\right).
   \]
   Thus choosing suitable $g$ we may assume $A = I_n$ and we have the normal form for $p_3$.
   
   Note that if $g$ is an orthogonal matrix $g^T  g = I_n$, then the action of 
   $
    \begin{bmatrix}
        g^{-1} &    0\\
        0      & g^T
    \end{bmatrix}
   $
   preserves $p_1$, $p_2$ and $p_3$.
   Thus it remains to prove that for $A$ a general symmetric matrix,
   there exists an orthogonal matrix $g$ such that $g^T A g$ is diagonal, 
   or in other words that two general quadratic polynomials can be simultaneously diagonalised.
   This is a standard fact, see for instance \cite[Prop.~2.1(a)\&(d)]{miles}.
\end{prf}

\section{Embedding and parametrisation}

There is a canonical morphism $\Wedge{n-2} V \stackrel{\wedge \omega}{\longrightarrow} \Wedge{n}V$
(taking into account that $\omega \in \Wedge{2} V^*$ determines a natural isomorphism between $V$ and $V^*$).
In \cite[Chapter~9, \S5, n.3]{bourbaki} it is shown that this morphism is injective,
and that there exists a canonical direct summand of the image which is exactly the weight space $V_{\omega_n}$.
(Here $V_{\omega_n}$ denotes the irreducible representation of $\Sp(2n)$ with highest weight $\omega_n$, 
the fundamental weight associated to the last simple root $\alpha_n$, see \cite{bourbakitables} for details). 
 In other words we have a splitting:
$$
  \Wedge{n} V=V_{\omega_n}\oplus \Wedge{n-2} V.
$$
The image of the Pl\"u{}cker embedding of the Lagrangian Grassmannian $\LG(n,2n)$ spans $\PP(V_{\omega_n})$ 
and we have the diagram:
\[
  \xymatrix{  & \PP(\Wedge{n}\CC^{2n})\\
  \LG(n,2n) \ar@{^{(}->}[ur] \ar@{^{(}->}[r]^{\psi} &\PP(V_{\omega_n}) \ar@{^{(}->}[u]}
\]
In the neighborhood of 
$\theta \left( \begin{bmatrix}
                I_n\\
		0
	\end{bmatrix} \right)
$ 
the embedding $\psi\colon\LG(n,2n) \hookrightarrow \PP(\Wedge{n}\CC^{2n})$ is given by sending 
$
\theta \left( \begin{bmatrix}
                I_n\\
		A
	\end{bmatrix} \right)
$ 
to all the possible $k \times k$ minors of the 
$n \times n$ symmetric matrix $A$:
\begin{align}
\label{map psi} \psi \colon \LG(n,2n) &\hookrightarrow \PP(V_{\omega_n}) \subset \PP(\Wedge{n}\CC^{2n})\\
\nonumber \theta \left(
	  \begin{bmatrix}
                I_n\\
		A
          \end{bmatrix} \right) &\mapsto [\:A_{IJ}\:]
\end{align}
with $I,J \subset \setfromto{1}{n}$ and $|I|=|J|=k$, with the convention that the $0 \times 0$ minor is just equal to $1$.
Notice also that $A_{IJ}=A_{JI}$.

It is useful to order the coordinates of $\PP(\Wedge{n}\CC^{2n})$ in an increasing order, so the index $k$ runs from $0$ to $n$:
$$[1,A,\Wedge{2} A, \ldots, \Wedge{n-1} A, \det A].$$
In this order the analogous neighborhood of 
$
 \theta \left(
\begin{bmatrix}
          0\\
		I_n
	\end{bmatrix}
\right)$ 
consisting of the classes of points
$
 \theta \left(
\begin{bmatrix}
         A\\
		I_n
	\end{bmatrix}\right)
$
is described in a symmetric way:
$$[\det A, \Wedge{n-1} A^T, \ldots, A^T,1]$$
with the appropriate choices of order and signs of minors. To see this, consider a point which is in both neighborhoods 
$
 \theta \left(\begin{bmatrix}
         I_n\\
		A
	\end{bmatrix} \right)= 
 \theta \left(	\begin{bmatrix}
	        A^{-1}\\
			I_n
		\end{bmatrix}\right)
$
with $A$ invertible. Then this point is mapped to:
\begin{align}
	\nonumber \big{[}1,A,\Wedge{2} A, \ldots, &\Wedge{n-1}A, \det A\big{]}=\\
	\nonumber &=\Big{[}\frac{1}{\det A},\frac{1}{\det A}A,\frac{1}{\det A}\Wedge{2} A,\ldots, \frac{1}{\det A}\Wedge{n-1}A,1 \Big{]}\\
	\nonumber &=\Big{[}\det (A^T)^{-1}, \Wedge{n-1}((A^T)^{-1}), \Wedge{n-2}((A^T)^{-1}),\ldots,(A^T)^{-1},1\Big{]}.
\end{align}

Here the equality $\frac{1}{\det A} \Wedge{k} A = \Wedge{n-k}((A^T)^{-1})$ is standard and well known, 
but rarely explicitly written down. See \cite[Prop.~H.19]{jabu_dr}, where $A$ is assumed to have determinant one 
(in the proof sketched there one can easily take into account an arbitrary determinant).

In $\Wedge{n} V = \Wedge{n} (W \oplus W^*)$ we distinguish the symmetric part:
\begin{align*}
  \Wedge{n}(W \oplus W^*)&=\bigoplus_{i=0}^n (\Wedge{i} W \otimes \Wedge{n-i}W^*)\\
                         &=\bigoplus_{i=0}^n (\Wedge{i} W \otimes \Wedge{i}W) \otimes \Wedge{n} W\\
                         &=\bigoplus_{i=0}^n (S^2(\Wedge{i}W) \oplus \Wedge{2}(\Wedge{i} W)) \otimes \Wedge{n} W.
\end{align*}
Since we are interested in the projectivisation of this vector space, the twist by $\Wedge{n} W$
becomes irrelevant and we regularly skip it.
The space $V_{\omega_n}$ is always contained in the symmetric part
$\bigoplus_{i=0}^n S^2(\Wedge{i}W)$, however (for $n\ge 4$) it is strictly smaller.

For simpler notation, we will consider all the symmetric minors, rather than its subset.
By this we mean that rather than working with  the embedding \eqref{map psi}:
$$
\psi:\LG(n,2n) \hookrightarrow \PP(V_{\omega_n}),
$$ 
we work with embedding $\varphi$:
\begin{equation}\label{map varphi}
	\varphi: \LG(n,2n)\hookrightarrow \PP(\bigoplus_{i=0}^n (S^2(\Wedge{i}W)).
\end{equation}

\section{Tangent space calculation}

Given $I,J \subseteq \setfromto{1}{n}$ with $|I|=|J|=k$ and given $A_{IJ}$ the corresponding $k\times k$ minor of 
the $n \times n$ symmetric matrix $A$, for any chosen value of $k$, define:
$$A_{IJ}^{ij}:=\Big[(-1)^{\#(i,I)+\#(j,J)}A_{I\setminus \{i\} J\setminus \{j\}} 
                      + (-1)^{\#(i,J)+\#(j,I)}A_{I\setminus \{j\} J\setminus \{i\}}\Big],$$
with the convention that $A_{I\setminus \{i\} J\setminus \{j\}}=0$ whenever 
$i \notin I$ or $j \notin J$ 
and where by $\#(i,I)$ we denote the index of $i$ in $I$, that is $\#\set{i' \in I \mid i'\le i}$.

\begin{lem}\label{jacobian 1pt}
The affine tangent space  $\hat{\Tang}_p(\LG(n,2n))$ to $\LG(n,2n)$ at the point 
$p=\theta\left(
   \begin{bmatrix}
                I_n\\
		A
	\end{bmatrix}\right)$
as a subspace of $\bigoplus_{i=0}^n \left(S^2(\Wedge{i}W)\right)$ is generated by the $d+1$ rows of the matrix:
$$
\left[
\begin{array}{c|c|c}
1&\:\:A\:\:&\Wedge{2} A, \Wedge{3} A,\:\ldots \:,\Wedge{n-1}A,\det A\\
\hline
0&&\\
\vdots&\:\:\:\:\:I_d\:\:\:\:\:&A_{IJ}^{ij}\\
0&&
\end{array}
\right]
$$
The rows of this matrix are indexed by unordered pairs $(i,j)$, with $1\le i,j \le n$,
and one extra row on the very top of the matrix.
The columns are indexed by unordered pairs $I,J$ of subsets of $\setfromto{1}{n}$, with $|I| =|J|\in 
\setfromto{2}{n}$.
\end{lem}

\begin{proof}
The affine tangent space $\hat{\Tang}_p \LG(n,2n)$ at the point 
$p=\theta\left(
   \begin{bmatrix}
                I_n\\
		A
	\end{bmatrix}\right)$
as a subspace of $\bigoplus_{i=0}^n \left(S^2(\Wedge{i}W)\right)$ is generated by
$\varphi(p)$ and the partial derivatives $(\pder{a_{ij}} \varphi)|_p$. We have:

\begin{align}
\nonumber &\left(\pder{a_{ii}} \varphi|_p\right)_{II}=\left\{
\begin{array}{ll}
0& i \notin I\\
\star& i \in I
\end{array}
\right.\\
\nonumber &\left(\pder{a_{ii}} \varphi|_p\right)_{IJ}=\left\{
\begin{array}{ll}
0& i\notin I\\
0& i \notin J\\
\star& i \in I \cap J
\end{array}
\right.\\
\nonumber &\left(\pder{a_{ij}} \varphi|_p\right)_{II}=\left\{
\begin{array}{ll}
0& i\notin I\\
0& j \notin I\\
\star& i, j \in I
\end{array}
\right.\\
\nonumber &\left(\pder{a_{ij}} \varphi|_p\right)_{IJ}=\left\{
\begin{array}{ll}
0& i\notin I\\
0& j \notin J\\
\star& i\in I \text { and } j \in J
\end{array}
\right.
\end{align}
where the symbol $\star$ is just a placeholder for a non-zero derivative.
By expanding the determinant either by row or column we explicitely compute the derivatives $\star$: 
$$\left(\pder{a_{ij}} \varphi|_p\right)_{IJ}=A_{IJ}^{ij}.$$
\end{proof}

\begin{lem}\label{jacobian diagonal}
In the notation of Lemma~\ref{jacobian 1pt},
 if the symmetric matrix $A$ is a diagonal matrix
 then $A_{IJ}^{ij}=0$, unless $I = K \cup \{i\}$ and $J = K \cup \{j\}$, 
 for some subset $K \subset \setfromto{1}{n}$.
\end{lem}

\begin{proof}
It is an immediate application of the above computations.
\end{proof}


\begin{lem}\label{tg space p1 and p2}
	In the proposed coordinates, for the points $p_1$ and $p_2$ as in Proposition~\ref{prop normal form} we have the following equalities 
	for affine tangent spaces as subspaces of $\bigoplus_{i=0}^n (S^2(\Wedge{i}W))$:
	\begin{align}
		\nonumber &\hat{\Tang}_{p_1}\LG(n,2n)=S^2(\Wedge{0} W) \oplus S^2(\Wedge{1} W) \simeq \CC \oplus S^2W\\
		\nonumber &\hat{\Tang}_{p_2}\LG(n,2n)=S^2(\Wedge{n-1} W) \oplus S^2(\Wedge{n} W) \simeq (S^2W^* \oplus \CC) \otimes (\Wedge{n} W)^{\otimes 2}
	\end{align}
	In particular, by Terracini Lemma~\ref{Terracini} 
        the second secant variety $\sigma_2(\LG(n,2n))$ always has the expected dimension. 
\end{lem}

\begin{prf}
	Immediate from the given parametrisations around the points $p_1$ and $p_2$.
\end{prf}

Slightly more demanding is the computation of the tangent space at the points $p_3$ and $p_4$,
 still from Proposition~\ref{prop normal form}. 

In the proposition below, for fixed $k$ we divide the $k \times k$ minors of $A$ into 3 groups. 
First the ``on-diagonal'' ones, namely those of the form $A_{II}$.
Then ``slightly off-diagonal'', namely those of the form $A_{IJ}$,
where $I =K \cup \set{l}$ and $J = K \cup \set{m}$, for some subset $K \subset \setfromto{1}{n}$,
with $l, m \in \setfromto{1}{n} \setminus K$ and $l \ne m$.
Finally, all the other minors, which in our setup (the choice of points in their normal form),
are irrelevant.

\begin{prop}\label{3 pts lemma}
   Let $p_1$, $p_2$, $p_3$  and $p_4$ be four general points of $\LG(n, 2n)$
   in their normal forms as in Proposition~\ref{prop normal form}.
   Then the space 
   \[
      \hat{\Tang}_{p_1}(\LG(n,2n)) + \hat{\Tang}_{p_2}(\LG(n,2n)) + \hat{\Tang}_{p_3}(\LG(n,2n)) + \hat{\Tang}_{p_4}(\LG(n,2n))
   \]
   is spanned by the rows of the matrix:
   \begin{equation}\label{3 pts jac matrix}
      \begin{array}{|c|c|c|c|c|c|c|c|}
        \hline
         1 & 0 \dotsc 0 & 0 \dotsc 0 & \dotsc\dotsc & 0 \dotsc 0& 0\dotsc 0 & 0\\
        \hline
&&&&&& \\
        \underline{0}&\:\:I_d\:\:&\underline{0}&\dotsc\dotsc &\underline{0}&\underline{0}&\underline{0}\\
&&&&&&\\
\hline
0&0 \ldots 0 & 0 \dotsc 0 & \dotsc\dotsc & 0 \dotsc 0 &0\ldots 0&1\\
        \hline
&&&&&&\\
        \underline{0}&\underline{0}&\underline{0}&\dotsc\dotsc &\underline{0}&\:\:I_d\:\:&\underline{0}\\
&&&&&&\\
        \hline
1& I_n & \Wedge{2}I_n & \dotsc\dotsc & \Wedge{n-2}I_n & \Wedge{n-1}I_n&1\\
\hline
&&&&&&\\
        \ast&\ast&  M_2 & \dotsc\dotsc & M_{n-2}
&\ast&\ast\\
&\phantom{\Wedge{n-1}Q_n}&\quad \phantom{\Wedge{n-2}Q_n} \quad&&&&\\
        \hline
1& Q_n & \Wedge{2}Q_n & \dotsc\dotsc & \quad \Wedge{n-2}Q_n \quad & \Wedge{n-1}Q_n& \det Q_n\\
\hline
&&&&&&\\
        \ast&\ast&  N_2 & \dotsc\dotsc & N_{n-2}
&\ast&\ast\\
&&&&&&\\
        \hline
      \end{array}
\end{equation}
where the $*$ are some matrices of the appropriate size,
and the matrices $M_k$ and $N_k$ consist of the following blocks:
\begin{align}\label{submatrix M}
M_k\:\: &:= \:\:\begin{array}{r||c|c|c|}
             &A_{II}            &A_{K \cup \{\ell\}K \cup \{m\}}         & A_{IJ}\\
\hline\hline
&&&\\
\pder{a_{ii}}|_{p_3}           &\:\:\:\:\:\:\:\:\:\:\:\:\:\:\delta_I^i\:\:\:\:\:\:\:\:\:\:\:\:\:\: & \underline{0}   &\underline{0}\\
&&&\\
\hline
&&&\\
\pder{a_{ij}}|_{p_3}    &\underline{0}    &\:\:\:\:\:\:\:\:\:\:\: (-1)^\epsilon \delta_{i\ell}\delta_{jm}  \:\:\:\:\: \:\: \:\:\:         & \underline{0}\\
&&&\\
\hline
\end{array}\\
\nonumber\\
\label{submatrix N}
N_k\:\: &:=\:\:\begin{array}{r||c|c|c|}
             &A_{II}            &A_{K \cup \{\ell\}K \cup \{m\}}         & A_{IJ}\\
\hline\hline
&&&\\
\pder{a_{ii}}|_{p_4}         &\:\:\: \prod_{j \in I\setminus \{i\}} q_j\delta_I^i\:\:\: & \underline{0}                                   &\underline{0}\\
&&&\\
\hline
&&&\\
\pder{a_{ij}}|_{p_4}    &\underline{0}    &(-1)^{\epsilon}\delta_{i\ell}\delta_{jm}  (\prod_{\beta \in K} q_\beta)              & \underline{0}\\
&&&\\
\hline
\end{array}
\end{align} 
The index $k$ ranges from $2$ to $n-2$, $|I|=|J|=k$, $|I\cap J| \leq k-2$.
The symbol $\delta_{ij}$ is $1$ if $i=j$ and $0$ otherwise. 
The symbol $\delta^i_I$ is $1$ if $i \in I$ and $0$ otherwise. 
Moreover by symmetry we can assume $\ell < m$ and $i<j$. 
Here $\epsilon = \#(i,K\cup \{\ell\})+\#(j,K\cup \{m\})$.
\end{prop}

\begin{proof} The proposition follows as a corollary of Lemmas~\ref{jacobian 1pt},
 \ref{jacobian diagonal} and \ref{tg space p1 and p2}. 

From Lemma~\ref{tg space p1 and p2} it is clear that the first $2d+2$ rows of the stacked Jacobian matrix will have the form above. 
We also conclude that from having our first two points $p_1$ and $p_2$ in their normal forms, we do not care about the first and last columns, 
in the lower part of the stacked matrix, hence the notation $*$.

Now let us look at the lower part of the stacked Jacobian matrix. We have that 
$p_3=
\theta\left( 
\begin{bmatrix}
         I_n\\
I_n
	\end{bmatrix}\right)
$
 and 
$p_4=
\theta\left( 
\begin{bmatrix}
         I_n\\
Q_n
	\end{bmatrix}\right)
$
with a diagonal matrix $Q_n$, 
so we can apply Lemma~\ref{jacobian diagonal} to both these points.
If the matrix $A$ is diagonal then in particular the minors 
$A_{IJ}^{ij}$ are nonzero if and only if either $I=J$, $i=j$ and $i \in I$ 
or else $I=K \cup \{i\}$ and $J=K\cup \{j\}$.
The explicit application of Lemma~\ref{jacobian 1pt} for
$A=I_n$ and $A=Q_n=\hbox{diag}(q_{11},\ldots,q_{nn})$ concludes the proof.
\end{proof}

Carefully looking at matrices \eqref{3 pts jac matrix}, \eqref{submatrix M} and \eqref{submatrix N}
 we conclude that in order to calculate the rank of \eqref{3 pts jac matrix}
 (equivalently, the dimension of the fourth secant variety) the following holds:
\begin{enumerate}
 \item we may use the first $2(d+1)$ rows to eliminate the $\ast$ parts of the lower rows, and
       thus the rank in question is equal to $2(d+1)$ plus the rank of the submatrix of the matrix 
       \eqref{3 pts jac matrix} obtained by removing the first $2(d+1)$ rows and the first and last $(d+1)$ columns. 
       We call this submatrix $B$.
 \item The submatrix $B$ is a direct sum of two matrices, $B^{diag}$, consisting of the rows corresponding 
       to $p_3, \pder{a_{ii}}|_{p_3}, p_4, \pder{a_{ii}}|_{p_4}$ 
       and the columns corresponding to minors $A_{II}$, and $B^{off}$, 
       consisting of the remaining rows and columns.
       Thus 
       \[
         \rk B = \rk B^{diag}+ \rk B^{off}.
       \]
\end{enumerate}


\section{Third secant variety}\label{section third secant case}


From Lemma~\ref{tg space p1 and p2} we know that the second secant variety $\sigma_2(\LG(n,2n))$
 always has the expected dimension. For $n \le 3$ this second secant fills the ambient space.
Thus for the rest of the paper we assume that $n \ge 4$.
In this section we calculate dimensions of third secant varieties.  

\begin{thm}
Suppose $n \ge 4$, and $r =3$.
Then:
\begin{itemize}
 \item If $n=4$, then $\sigma_3(\LG(4,8))$ has defect 1.
 \item If $n\ge 5$, then $\sigma_3(\LG(n,2n))$ has the expected dimension.
\end{itemize}
\end{thm}

\begin{proof}
By Terracini Lemma~\ref{Terracini} and Proposition~\ref{3 pts lemma},
  we need to calculate the rank of the first three blocks of the stacked Jacobian matrix \eqref{3 pts jac matrix}, 
  i.e., the blocks corresponding to the points $p_1$, $p_2$ and $p_3$ and the respective derivatives 
  (so the first $3d+3$ rows).
The expected rank is the maximal one, $3d+3$, which is the expected dimension of 
the affine cone of $\sigma_3(\LG(n,2n))$.

The first $2d+2$ rows are linearly independent, so we focus our attention on the third block of $d+1$ rows,
and in particular on the submatrices $M_k$ described in \eqref{submatrix M}.
We need to show that the rows, restricted to the columns corresponding to $k =\fromto{2}{n-2}$
have maximal rank $d+1$ for $n \ge 5$ and rank $d = 10$ for $n=4$. 
In case $n=4$ the index $k$ can only be equal to $2$, whereas for $n\ge 5$,
we at least have the blocks $k=2$ and $k=3$ at our disposal.
We will show that the rank of the block corresponding to $k=2$ is equal to $d$,
but the unique (up to scale)
linear relation between the rows does not extend to the block corresponding to $k=3$.

We ``zoom in'' the two blocks $M_2$ and $M_3$ of the matrix described in \eqref{submatrix M}, 
together with the row coming from the point $p_3$.
Recall once again that we are ordering the minors by putting first the on-diagonal $A_{II}$ 
ones and then the off-diagonal $A_{IJ}$.
We also order the elements of the symmetric matrix $A$ in the same on-diagonal, 
off-diagonal order, so that the on-diagonal block has $n$ rows, 
and the off-diagonal block has the other $d-n$ rows.

\begin{equation}\label{1st zoom in}
\begin{array}{r||c|c|c||c|c|c|}
&\multicolumn{2}{r}{ M_2\phantom{blabl}}&&\multicolumn{2}{r}{ M_3\phantom{blabl}}&\\
&
A_{\{\ell m\}\{\ell m\}} &A_{\{\beta\ell\}\{\beta m\}} & A_{\{\beta\ell\}\{\gamma m\}} &A_{II} &A_{K \cup \{\ell\}K \cup \{m\}} & A_{IJ}\\
&&&&&|K|=2&\\
\hline\hline
1&1 \ldots 1&0 \ldots \ldots 0&0 \ldots 0&1\ldots 1&0 \ldots \ldots 0 &0 \ldots 0\\
\hline
&&&&&&\\
\pder{a_{ii}} &\:\:\:\delta^i_{\{\ell m\}}\:\:\: & \underline{0} &\underline{0}&\:\:\delta_I^i\:\:&\underline{0}&\underline{0}\\
&&&&&&\\
\hline
&&&&&&\\
\pder{a_{ij}}&\underline{0} &(-1)^\epsilon \delta_{i\ell}\delta_{jm} & \underline{0}&\underline{0}&(-1)^\epsilon \delta_{i\ell}\delta_{jm}&\underline{0}\\
&&&&&&\\
\hline
\end{array}
\end{equation}

In the first place,  let us focus only on $M_2$.
The off-diagonal part is of maximal rank $d-n$ because there is precisely one non-zero entry in every column
and every row is non-zero. 
On the other hand the sum of the rows of the on-diagonal block equals twice the first row, 
the one corresponding to the point $p_3$, so we have a linear relation.
We claim that the on-diagonal block of $M_2$ (without the first row corresponding to $p_3$)
 is of maximal rank $n$.
Let's ``zoom-in'' the block $M_2$ even more, focusing on the on-diagonal part:

\begin{equation}\label{2nd zoom in}
\begin{array}{r||c|c|c}
&A_{\{1\ell\}\{1\ell\}}&A_{\{2\ell\}\{2\ell\}}&\:\:\:\ldots\:\:\: \\
&\ell=2,\ldots,n&\ell=3,\ldots,n&\ldots\\
\hline
\hline       
\pder{a_{11}}&1\ldots \ldots 1&&\\
\hline
\pder{a_{22}}&&1\ldots \ldots 1&\\
\vdots&I_{n-1}&&\ddots\\  
\vdots&&I_{n-2}&\ddots\\  
\pder{a_{nn}}&&& 
\end{array}
\end{equation}
Notice there is a copy of $I_{n-1}$, so the rank of the diagonal block in $M_2$ is at least $n-1$.
But if it is of rank $n-1$, then the first row of $M_2$ (the one corresponding to $\pder{a_{11}}$)
is the sum of all the other rows. One easily verifies this is not the case for instance 
on the columns corresponding to $A_{\{2\ell\}\{2\ell\}}$.
Thus $\rk M_2 = n$ and this finishes the proof of the $n=4$ case. 

Suppose $n \ge 5$, so that for $k=3$ we have $k\le n-2$.
By the above argument, the rank of the matrix \eqref{1st zoom in} is at least $d$. 
If it is equal to $d$,
then the aforementioned relation between rows of $M_2$
and the first row holds also for the part of the matrix corresponding to the $3 \times 3$ minors.
This is not the case for any column corresponding to $A_{II}$ ---
in such a column there are exactly $3$ values of $i$, where $\delta^i_I$ is $1$, otherwise $\delta^i_I$ is $0$.
Thus the rank of the matrix \eqref{1st zoom in} is $d+1$ and $\sigma_3(\LG(n,2n))$ is non-defective.
\end{proof}


\section{Fourth secant variety}\label{section fourth secant case}


Continuing the proof of Theorem~\ref{thm third and fourth secant}
 we focus on the remaining case $r=4$.
 
\begin{thm}
Suppose $n \ge 4$, and $r =4$.
Then:
\begin{itemize}
 \item If $n=4$, then $\dim \sigma_4(\LG(4,8))$ has defect 4.
 \item If $n\ge 5$, then $\sigma_4(\LG(n,2n))$ has the expected dimension.
\end{itemize}
\end{thm}

\begin{proof}
We use the fourth point $p_4$ in the normal form. 
This time we need to show that the rank of the whole stacked Jacobian
matrix \eqref{3 pts jac matrix} is equal to $40$ for $n=4$ and $4d+4$ for $n \ge 5$.

Similarly to the case of $r=3$, we focus our attention on the last $2d+2$ rows,
and on the middle columns, i.e.~on the $M_k$'s and $N_k$'s for $k \in \setfromto{2}{n-2}$.
We first claim that $2(d-n)$ rows corresponding to the off-diagonal $\pder{a_{ij}}$ 
at both points $p_3$ and $p_4$ are always linearly independent.
Indeed, it is enough to look at the appropriate parts of $M_2$ and $N_2$.
Fixing a column $A_{\set{\beta \ell} \set{\beta m}}$ there are precisely $2$ non-zero entries,
the ones corresponding to $\pder{a_{\ell m}}$ on both points.
Thus the only possible way that there is a linear relation between the $2(d-n)$ rows 
is that $\pder{a_{\ell m}}|_{p_3}$ is proportional to $\pder{a_{\ell m}}|_{p_4}$.
We write explicitly the (non-zero columns of) two rows to observe that this is not the case.
We ignore the sign, as within a column it is the same sign, thus we may multiply the column by $-1$ 
if necessary and this obviously does not change the rank.
\[
\begin{array}{r||c|c|c|c}
&A_{\set{1\ell}\set{1 m}}& \ldots & A_{\set{\beta\ell }\set{\beta m}}& \ldots \\
\hline
\hline
\pder{a_{\ell m}}|_{p_3}& 1        & \ldots &               1  & \ldots \\
\pder{a_{\ell m}}|_{p_4}& q_1   && q_\beta  &  \\
\end{array}
\]
Here $\beta$ runs through $\setfromto{1}{n}\setminus \set{\ell, m}$.
Clearly, if $\ell$ or $m$ is $1$, then the first column does not show up, but anyway,
since $n \ge 4$, there are at least two possible values of $\beta$
and as long as $q_i$ are not pairwise equal, then the two rows are independent.
Since the $q_i$ are general, indeed the $2(d-n)$ off-diagonal rows are independent.

It remains to look at the $2n$ diagonal rows and the $2$ rows corresponding to $p_3$ and $p_4$.

First suppose $n=4$. Then we are able to write the matrix $B^{diag}$ explicitly:
\[
  \begin{pmatrix}
    1  & 1  & 1  & 1  & 1  &  1  \\
    1  & 1  & 1  & 0  & 0  &  0  \\
    1  & 0  & 0  & 1  & 1  &  0  \\
    0  & 1  & 0  & 1  & 0  &  1  \\
    0  & 0  & 1  & 0  & 1  &  1  \\
q_1 q_2 & q_1 q_3  & q_1 q_4 &q_2 q_3 &q_2 q_4 &q_3 q_4 \\
       q_2 &        q_3  &        q_4 &          0   &           0  &            0 \\
q_1        &             0  &             0 &       q_3 &       q_4 &            0 \\
            0 & q_1         &             0 &q_2        &            0 &       q_4 \\
            0 &             0  & q_1        &            0 &q_2        &q_3        \\
    
  \end{pmatrix}.
\]
Since the matrix has only $6$ columns (and $10$ rows) it is impossible that the rows are linearly independent.
Using Gaussian elimination one can conclude that this matrix has in fact rank $6$ (for mildly general $q_i$'s),
 and thus the dimension of the secant variety is $2(d+1) + 2(d-n) + 6 - 1 = 2(10+1) + 2(10 -4) + 5 = 39$.

For the case $n \geq 5$ we want to prove that the rank of the submatrix $B^{diag}$ is exactly $2n$
(let us ignore for a second the $2$ rows corresponding to $p_3$ and $p_4$, and simply look at the $2n$ diagonal rows).
$B^{diag}$ will consist of a copy of the submatrix (\ref{2nd zoom in})
that we have already described, stacked above a matrix of the form:

\begin{equation}\label{2nd zoom in for the diagonal}
\begin{array}{r||c|c|c}
&A_{\{1\ell\}\{1\ell\}}&A_{\{2\ell\}\{2\ell\}}&\:\:\:\ldots\:\:\: \\
&\ell=2,\ldots,n&\ell=3,\ldots,n&\ldots\\
\hline
\hline       
\pder{a_{11}}&q_2\:q_3\:q_4\ldots \ldots q_n&&\\
\hline
\pder{a_{22}}&&q_3\:q_4\ldots \ldots q_n&\\
\vdots&q_1I_{n-1}&&\ddots\\  
\vdots&&q_2I_{n-2}&\ddots\\
&&&\\ 
\pder{a_{nn}}&&& 
\end{array}
\end{equation}

Let us change the order of the rows so that the corresponding derivatives of (\ref{2nd zoom in}) and (\ref{2nd zoom in for the diagonal}) are 
adjacent. For clarity, we rewrite the matrix $B^{diag}$ here below.
{\footnotesize
$$
\begin{array}{|cccc|cccc|cccc|c|cc|c|}
*      &*   &       &    &*   &*   &       &    &*   &*   &       &    &\ldots       &        &        &\\
	\hline
1      &1   &\ldots &1   &    &    &       &    &    &    &       &    &             &        &        &\\
q_2    &q_3 &\ldots &q_n &    &    &       &    &    &    &       &    &             &        &        &\\
1      &    &       &    &1   &1   &\ldots &1   &    &    &       &    &             &        &        &\\
q_1    &    &       &    &q_3 &q_4 &\ldots &q_n &    &    &       &    &             &        &        &\\
       &1   &       &    &1   &    &       &    &1   &1   &\ldots &1   &             &        &        &\\
       &q_1 &       &    &q_2 &    &       &    &q_4 &q_5 &\ldots &q_n &             &        &        &\\
       &    &       &    &    &1   &       &    &1   &    &       &    &\ldots\ldots &        &        &\\
       &    &       &    &    &q_2 &       &    &q_3 &    &       &    &\ldots\ldots &        &        &\\
       &    &\ddots &    &    &    &\ddots &    &    &1   &       &    &             &1       &1       &\\
       &    &\ddots &    &    &    &\ddots &    &    &q_3 &       &    &\ddots       &q_{n-1} &q_n     &\\
       &    &       &    &    &    &       &    &    &    &\ddots &    &\ddots       &1       &        &1\\
       &    &       &    &    &    &       &    &    &    &\ddots &    &             &q_{n-2} &        &q_n\\
       &    &       &1   &    &    &       &1   &    &    &       &1   &\ldots\ldots &        &1       &1\\
       &    &       &q_1 &    &    &       &q_2 &    &    &       &q_3 &\ldots\ldots &        &q_{n-2} &q_{n-1}\\
\hline
\end{array}
$$
}
We are now going to focus our attention on the first $2$ columns of the first $n-5$ blocks,
namely the columns:
\[
  \fromto{A_{\set{1,2}\set{1,2}}, A_{\set{1,3}\set{1,3}}, A_{\set{2,3}\set{2,3}},A_{\set{2,4}\set{2,4}}}
   {A_{\set{n-5,n-3}\set{n-5, n-3}}}. 
\]
Some of those columns are marked by $*$ above. 
We also consider the last 10 columns. 
Altogether these columns have a following lower triangular block form:
{\footnotesize
\[
\begin{array}{|cc|cc|cc|cc|ccc|}
	\hline
1      &1      &       &       &       &       &        &        &    &         &\\
q_2    &q_3    &       &       &       &       &        &        &    &         &\\
\vdots &\vdots &1      &1      &       &       &        &        &    &         &\\
       &       &q_3    &q_4    &       &       &        &        &    &         &\\
       &       &\vdots &\vdots &\ddots &       &        &        &    &         &\\
       &       &       &       &\ddots &\ddots &        &        &    &         &\\
       &       &       &       &\vdots &\vdots &1       &1       &    &         &\\
       &       &       &       &       &       &q_{n-4} &q_{n-3} &    &         &\\
\hline
       &       &       &       &       &       &\vdots  &\vdots  &    &         & \\
       &       &       &       &       &       &        &        &    &\Upsilon & \\
       &       &       &       &       &       &        &        &    &         & \\
\hline
\end{array}
\]
}
where $\Upsilon$ is the following $10\times 10$ matrix:
\[
\begin{array}{|cccc|ccc|cc|c|}
\RotateUp{\set{n-4,n-3}} &\RotateUp{\{n-4,n-2\}} &
\RotateUp{\set{n-4,n-1}} &\RotateUp{\{n-4,n\}} &\RotateUp{\{n-3,n-2\}} &\RotateUp{\{n-3,n-1\}} & 
      \RotateUp{\{n-3,n\}}  &\RotateUp{\{n-2,n-1\}}  &\RotateUp{\{n-2,n\}} &\RotateUp{\{n-1,n\}}\\
\hline  
1&1&1&1&0&0&0&0&0&0 \\
q_{n-3}&q_{n-2}&q_{n-1}&q_n&0&0&0&0&0&0 \\
1&0&0 & 0&1&1&1      &0       &0      &0\\
q_{n-4}&0&0 & 0&q_{n-2}&q_{n-1}&q_n      &0       &0      &0\\
0&1&0 & 0&1&0&0      &1       &1      &0\\
0&q_{n-4}&0 &0&q_{n-3}&0&0       &q_{n-1}  &q_n     &0\\
0&0&1 &0&0&1&0       &1       &0       &1\\
0&0&q_{n-4} &0&0&q_{n-3}&0       &q_{n-2}  &0       &q_n\\
0&0&0 &1&0&0&1       &0    &1       &1\\
0&0&0 &q_{n-4}&0&0&q_{n-3}      &0    &q_{n-2}  &q_{n-1}
\end{array}.
\]
Since the $q_i$'s are general, each of the blocks
$
  \begin{bmatrix}
     1  &1       \\
     q_i& q_{i+1}
  \end{bmatrix}
$
has rank $2$. 
Also it can be verified (either using a computer algebra system or by a patient Gaussian elimination) that
the rank of $\Upsilon$ is $10$. Thus $\rk B^{diag} = 2(n-5) + 10 = 2n$. 

So all in all the rank of the stacked matrix 
$\begin{bmatrix}
  M_2\\
  N_2
\end{bmatrix}$
is $2(d-n)$ from the off-diagonal part, plus $2n$ from the on-diagonal part, so in total $2(d-n)+2n=2d$.
So this means that there are 2 linear relations among the $2d+2$ rows. 
It is easy to see what they are (by analogy to the case $r=3$) and 
 these same relations (nor their linear combination) cannot hold on the block 
$\begin{bmatrix}
	M_3\\
	N_3
\end{bmatrix}$,
and this concludes the proof.
\end{proof}

\section{Computational experiments}\label{section_computations}

For small values of $n$ 
the secant varieties $\sigma_r(\LG(n,2n))$ all have the expected dimension,
except for the defective cases covered by Theorem~\ref{thm third and fourth secant}.

\begin{prop}\label{prop small values of n}
   Suppose $n \le \MaxVn$.
   Then $\sigma_r (\LG(n, 2n))$ have the expected dimension, 
   unless $n=4$ and $r \in \set{3,4}$.
\end{prop}
\begin{proof}
   The proof uses a naive computer code in Magma \cite{magma}.
   The code generates $r$ random points on $\LG(n, 2n)$ 
   (for a slight improvement of time needed to finish the calculation,
   the first $4$ points are assumed to be in the normal forms of Proposition~\ref{prop normal form}).
   Then it calculates the sum of the affine tangent spaces at these points.
   By Terracini Lemma~\ref{Terracini} and semicontinuity
   (since the rank can only drop at special points), 
   the dimension of the sum is a lower bound for the dimension of $\sigma_r(\LG(n,2n))$.

   For each $n \in \setfromto{4}{\MaxVn}$, we start the above experiment with $r=5$ 
   (Theorem~\ref{thm third and fourth secant} covers all the cases with $r \le 4$),
   and repeatedly add a new point until the dimension of the sum of tangent spaces
   is equal to the dimension of the ambient vector space.
   In all the cases we have obtained the lower bound to be equal to the expected dimension,
   which is an upper bound. Thus the dimension of the secant variety is the expected one for all these cases.
\end{proof}

\bibliographystyle{alpha}
\bibliography{references}

\end{document}